\documentclass[12pt]{amsart}
\usepackage{amsfonts}
\usepackage{amssymb}
\usepackage{latexsym}
\usepackage{setspace}
\usepackage{xcolor}

\newenvironment{thm}{\subsection{}{\textbf {Theorem.}}\em}{}
\newenvironment{prop}{\subsection{}{\textbf {Proposition.}}\em}{}
\newenvironment{cor}{\subsection{}{\textbf {Corollary.}}\em}{}
\newenvironment{lem}{\subsection{}{\textbf {Lemma.}}\em}{}

\newenvironment{pf}{\noindent{\em
{Proof.}}}{\begin{flushright}\eop \end{flushright}\smallskip}

\newenvironment{defn}{\subsection{}{\textbf 
{Definition.}}\em}{\smallskip}

\newenvironment{rem}{\subsection{}{\textbf {Remark.}}}{\smallskip}

\newcommand\cA{\ensuremath{\mathcal A}}

\newcommand\cJ{\ensuremath{\mathcal J}}

\newcommand\cM{\ensuremath{\mathcal M}}

\newcommand\cQ{\ensuremath{\mathcal Q}}

\newcommand\eop{{{\hfil \ensuremath \Box}}}
\begin{document}
\title{Invariant Subspaces in the  Dual of $A_{cb}(G)$ and $A_M(G)$}
\author{Brian Forrest, John Sawatzky and  Aasaimani Thamizhazhagan
 }

\maketitle

\begin{abstract} Let $G$ be a locally compact group. In this paper, we study various invariant subspaces of the duals of the algebras $A_M(G)$ and $A_{cb}(G)$ obtained by taking the closure of the Fourier algebra $A(G)$ in the multiplier algebra $MA(G)$ and completely bounded multiplier algebra $M_{cb}A(G)$ respectively. In particular, we will focus on various functorial properties and containment relationships between these various invariant subspaces including the space of uniformly continuous functionals and the almost periodic and weakly almost periodic functionals. 

Amongst other results, we show that if $\cA(G)$ is either $A_M(G)$ or $A_{cb}(G)$, 
then $UCB(\cA(G))\subseteq WAP(G)$ if and only if $G$ is discrete.  We also show that if $UCB(\cA(G))=\cA(G)^*$, then every amenable closed subgroup of $G$ is compact. 

Let $i:A(G)\to \cA(G)$ be the natural injection. We show that if $X$ is any closed topologically introverted subspace of $\cA(G)^*$ that contains $L^1(G)$, then $i^*(X)$ is closed in $A(G)$ if and only if $G$ is amenable.

\end{abstract}

\footnote{{\it Date}: \today.


This research is partially supported by the Natural Sciences and Engineering Research Council of Canada. }

\section{Introduction} Let $G$ be a locally compact group. The Fourier algebra of $G$ is an   algebra $A(G)$ of continuous functions on $G$ under pointwise operations consisting of the set of all coefficient functions $u(x)=\langle \lambda(x)f,g \rangle $ of the left regular representation of $G$ on the Hilbert space $L^2(G)$. $A(G)$ becomes a Banach algebra via the norm it inherits as the natural pre-dual of the group von Neumann algebra 
\[VN(G)=\overline{\textnormal{span}\{\lambda(x)|x\in G\}}^{-WOT}\subseteq B(L^2(G)) \]
where $\lambda(x)f(y)=f(x^{-1}y)$ for each $x,y\in G$ and $f\in L^2(G)$. 
If $G$ is abelian, then $A(G)$ is isometrically isomorphic as a Banach algebra with $L^1(\hat{G})$ via the classical Fourier transform and $VN(G)$ can be identified with $L^{\infty}(\hat{G})$, where $\hat{G}$ is the Pontryagin dual of $G$. 

For an abelian group $G$ there are several distinguished $C^*$-subalgebras of $L^{\infty}(\hat{G})$ that have been very well studied. These include the spaces of uniformly continuous functions $UCB(\hat{G})$, the almost periodic functions $AP(\hat{G})$, the weakly almost periodic functions $WAP(\hat{G})$ and $C_0(\hat{G})$, the space of continuous functions vanishing at infinity. All of these spaces induce corresponding subspaces of $VN(G)$ via the Fourier transform. In particular, we refer to the subspaces corresponding to  $UCB(\hat{G})$, $AP(\hat{G})$ and $WAP(\hat{G})$ in $VN(G)$ as the uniformly continuous, almost periodic and weakly almost periodic functionals on $A(G)$ respectively. The image of $C_0(\hat{G})$ is the reduced group $C^*$-algebra $C^*_{\lambda}(G)$. One common property that all of these spaces have in common is that they are invariant under the natural module action of $A(G)$ on $VN(G)$. 

When $G$ is non-abelian there is still a dual relationship between $A(G)$ and $L^1(G)$ but it is much more complex given the lack of a Pontryagin dual group. It is, however, possible to define analogous versions of all of the above subspaces of $VN(G)$ and indeed this can be done in the abstract for any commutative Banach algebra. For the Fourier algebra in particular, there is a significant body of work investigating the properties of the subspaces highlighted above with a focus on what their structures tell us about the underlying group. See for example \cite{ForMiao}, \cite{Gran}, \cite{KanLau}, \cite{Lau1} and \cite{Lau2}. 

In this paper, we will focus our attention primarily on non-amenable groups and on  analogous invariant subspaces in the duals of two Banach algebras that arise from the Fourier algebra, $A_M(G)$ and $A_{cb}(G)$, which are the closures of $A(G)$ in its multiplier and completely bounded multiplier algebra respectively. In particular, we will consider  various functorial properties and containment properties associated with these subspaces. Amongst other results,  if $\cA(G)$ is either $A_M(G)$ or $A_{cb}(G)$, we will characterize discrete groups as those locally compact groups for which $UCB(\cA(G))\subseteq WAP(\cA(G))$.  We also show that if $UCB(\cA(G))=\cA(G)^*$, then every amenable closed subgroup of $G$ is compact. In particular, $G$ has an open compact subgroup. 

$A(G)$ injects contractively into both $A_M(G)$ and $A_{cb}(G)$ with the range of the injection being proper precisely when $G$ is non-amenable. If  we let $i: A(G)\to \cA(G)$ be the natural injection, we will show that if $X$ is any of the closed subspaces  of $\cA(G)^*$ mentioned above  that contains $L^1(G)$, then $i^*(X)$ is closed in $A(G)$ if and only if $G$ is amenable.

\section{Preliminaries and Notation} 

Throughout this paper, $\cA$ will denote a  Banach algebra. In this case, the dual $\cA^*$ becomes a Banach $\cA$-bimodule with respect to the module actions 
\[\langle  u\cdot T, v \rangle =\langle T, vu \rangle ~~~~~\textnormal{and}~~~~~ 
\langle  T\Box u, v \rangle  = \langle T, uv \rangle  \]
for every $u,v\in \cA$ and $T \in \cA^*$. Note that when $\cA$ is commutative 
$\langle u\cdot T, v \rangle =\langle  T\Box u, v \rangle $ which will generally be the case for most of this paper. In fact, from this point we will assume that every Banach algebra $\cA$ we consider is commutative unless stated otherwise.

\begin{rem}
It is well known that there are two natural products that can be used to extend the multiplication of $\cA$ to its second dual $\cA^{**}$. In this paper we will choose the following Arens product: 
\begin{itemize}
\item[A1)] $\langle u\cdot T, v \rangle  =\langle T, vu \rangle $ for every $u,v \in \cA$ and $T \in \cA^*$.
\item[A2)] $\langle n\odot T, u \rangle  = \langle n, u \cdot T \rangle $ for every $u \in \cA$ and $T \in \cA^*$ and $n\in \cA^{**}$.
\item[A3)] $\langle m\odot n, T \rangle  = \langle m, n\odot T \rangle $ for every $T \in \cA^*$ and $m,n \in \cA^{**}$.
\end{itemize}

\end{rem}

 From here we will proceed with the following definitions and notational conventions. 

\begin{defn} We call the space 
\[UCB(\cA)=\overline{\textnormal{span}\{v \cdot T\mid v\in \mathcal{\cA}, T\in \cA^* \}}^{-\|\cdot\|_{\cA^*}}\]  
the \textit{uniformly continuous functionals on} $\cA$. 

We call $T\in \cA^*$ a \textit{(weakly) almost periodic functional} on $\cA$ if 
\[\{ u \cdot T\mid u \in \cA, \|u\|_{\cA} \leq 1\}\]
is relatively (weakly) compact in $\cA^*$ and we denote the space of all (weakly) almost periodic functionals on  $\cA$ by $AP(\cA)$ ($WAP(\cA)$).

\end{defn}

\begin{defn} We say that a closed subspace $X\subseteq \cA^*$ is \textit{invariant} if $u \cdot T \in X$ for every $u\in \cA$ and $T\in X$. \vspace{.05in}

Given a closed invariant  subspace $X$. We say that $X$ is topologically introverted if $m\odot T \in X$ for every $m\in X^*$ and $T \in X$.

\end{defn} 

\begin{rem} If $X$ is topologically introverted then $X^*$ can be made into a Banach algebra by mimicking what we did for $\cA^{**}$ by defining  
\[\langle m\odot n, T \rangle  = \langle m, n\odot T \rangle \]
for $m,n\in X^*$ and $T\in X$. 
    
It is a well-known criterion of Grothendieck that $T\in \cA^*$  is weakly almost periodic if and only if given two nets $\{u_{\alpha}\}_{\alpha \in \Omega_1}$ and  $\{v_{\beta}\}_{\beta\in \Omega_2}$ in $\cA$ we have that 
\[\lim\limits_{\alpha} \lim\limits_{\beta}\langle  T,u_{\alpha}v_{\beta}\rangle=
\lim\limits_{\beta} \lim\limits_{\alpha}\langle  T,u_{\alpha}v_{\beta}\rangle\]
whenever both limits exist. From this one can show that if $\cA$ is a commutative Banach algebra and $X$ is  topologically introverted, then $X^*$ is commutative if and only if $X\subseteq WAP(\cA)$. (See \cite{DuncHoss} and \cite{Pym}.)

\end{rem}

The following three propositions will prove useful.

\begin{prop} Let $\cA$ be a commutative Banach algebra. Then each of $UCB(\cA)$, $AP(\cA)$ and $WAP(\cA)$ are invariant subspaces of $\cA^*$.
\end{prop}

\begin{pf} That $UCB(\cA)$ is invariant is obvious since for each $T\in \cA^*$, we have $u\cdot T \in UCB(\cA)$. 

Let $T\in (W)AP(\cA)$. Let $v\in \cA$. Without loss 
of generality, we can assume that $\|v\|_{\cA}=1$. 
Since
\[\{ u \cdot(v\cdot T)\mid u \in \cA, \|u\|_{\cA} \leq 1\} \subset \{ u \cdot T\mid u \in \cA, \|u\|_{\cA} \leq 1\}\]
we have that $\{ u \cdot(v \cdot T)\mid u \in \cA, \|u\|_{\cA} \leq 1\}$ is relatively (weakly) compact and hence $v\cdot T\in (W)AP(\cA)$.
    
\end{pf}
\begin{prop} Let $\cA$ be a commutative Banach algebra. Then $UCB(\cA)$ is topologically introverted.
\end{prop}

\begin{pf} Let $m\in UCB(\cA)^*$ and let $T=u\cdot T_1\in UCB(\cA)$. Let $v\in \cA$. Then 
\begin{eqnarray*}
 \langle m\odot T,v \rangle &=&\langle m, v \cdot T \rangle   \\
 &=& \langle m, v\cdot (u\cdot T_1) \rangle  \\
 &=& \langle m,uv \cdot T_1 \rangle  \\
 &=& \langle m\odot T_1, uv \rangle \\
 &=& \langle u\cdot(m\odot T_1), v \rangle .\\
\end{eqnarray*}
It follows that $m\odot T= u\cdot(m\odot T_1)\in UCB(\cA)$. Moreover, if $T\in \textnormal{span}\{u\cdot T_1 | u\in \cA, T_1\in \cA^*\}$, then $m\odot T \in UCB(\cA)$ as well.

Finally, lets assume that $\{T_n\}\subset UCB(\cA)$ and that 
$\lim\limits_{n\to \infty}T_n =T $. If $u\in \cA$ with $\|u\|_{\cA}\leq 1$. Then for $m\in UCB(\cA)^*$, 
\[|\langle m\odot T, u \rangle -\langle m\odot T_n, u \rangle |=|\langle m, u\cdot (T-T_n) \rangle |\leq \|m\|_{UCB(\cA)^*}\|T-T_n\|_{\cA^*}. \] 
It follows that $\lim\limits_{n\to \infty} m\odot T_n =m\odot T $. Since $\textnormal{span}\{u\cdot T_1 | u\in \cA, T_1\in \cA^*\}$ is dense in $UCB(\cA)$, it follows that $m\odot T\in UCB(\cA)$ for every 
$m\in UCB(\cA)^*$ and $T\in UCB(\cA)$. That is, $UCB(\cA)$ is topologically introverted.

\end{pf}

\begin{prop} Let $X$ be a closed invariant subspace of $WAP(\cA)$. Then $X$ is topologically introverted. 
\end{prop}

\begin{pf} Let $X$ be a closed invariant subspace of $WAP(\cA)$. Let $T\in X$. Let 
\[O(T)=\{ u \cdot T\mid u \in \cA, \|u\|_{\cA} \leq 1\}.\]
Since $T\in WAP(\cA)$, the convex set $O(T)$ is 
relatively $\sigma(\cA^*,\cA^{**})$ compact. Hence its 
$\sigma(\cA^*,\cA^{**})$ closure $\overline{O(T)}^{\sigma(\cA^*,\cA^{**})}$ is weakly compact and hence is also a norm closed set in $X$. Moreover, the  $\sigma(\cA^*,\cA^{**})$ and $\sigma(\cA^*,\cA)$ topologies agree on $\overline{O(T)}^{\sigma(\cA^*,\cA^{**})}$. 

Next, we let $m\in X^*$. Without loss of generality, we may assume that $\|m\|_{X^*}=1$. Let $M\in \cA^{**}$ be an extension of $m$ to $\cA^*$ with 
$\|m\|_{X^*}=\|M\|_{\cA^*}=1$. Then by Goldstine's 
Theorem there exists a net $\{u_{\alpha}\}_{\alpha\in \Omega}$ in $\cA$ such that $\|u\|_{\alpha}\leq 1$ and $\{u_{\alpha}\}_{\alpha\in \Omega}$ converges to $M$ in the weak$^*$ topology on $\cA^{**}$. It follows that if $u\in\cA$ and $T\in X$, then 
\[\langle M\odot T, u \rangle = \langle M, u\cdot T \rangle =\lim\limits_{\alpha\in \Omega} \langle u_{\alpha}, u\cdot T \rangle =\lim\limits_{\alpha\in \Omega} \langle u_{\alpha}\cdot T, u \rangle .\] But this shows that 
  \[M\odot T\in \overline{O(T)}^{\sigma(\cA^*,\cA)}\subset X.\]  
  Moreover, if $u\in\cA$ and $T\in X$, then 
  \[\langle M\odot T, u \rangle = \langle M, u\cdot T \rangle =\langle m, u\cdot T \rangle =\langle m\odot T, u \rangle .\] Hence $m\odot T \in X$. It follows that $X$ is topologically introverted.
\end{pf}

\begin{defn} Let $I $ be a closed ideal in $\cA$. We let 
\[Z(I)=\{x\in \Delta(\cA)~|~u(x)=0~\textnormal{for all}~u \in I\}.\]

Given a closed set $E\subseteq \Delta(\cA)$, we let 
\[I(E)=\{u\in \cA~|~u(x)=0~\textnormal{for all}~x \in E\}.\]

We say that a closed set $E\subseteq \Delta(\cA)$ is a set of spectral synthesis for $\cA$ if the only closed ideal $I$ is $\cA$ with $Z(I)=E$ is $I(E)$.

\end{defn}

\section{Multipliers of the Fourier Algebra} 

Let $G$ be a locally compact group. We let $A(G)$ and $B(G)$ denote the Fourier and
 Fourier-Stieltjes algebras of $G$, which are Banach algebras of continuous functions
on $G$ and were introduced in [5]. 
 A \emph{multiplier} of $A(G)$ 
is a (necessarily bounded and continuous) function 
$v \!: G \to \mathbb{C}$
such that $v A(G) \subseteq A(G)$. For each multiplier $v$ of $A(G)$, 
the linear operator $M_v$ on $A(G)$ defined by 
$M_{v}(u)=vu$ for each $u\in A(G)$ is bounded via the Closed 
Graph Theorem. The {\it multiplier algebra\/} 
of $A(G)$ is the closed subalgebra

\[
  MA(G) := \{ M_v : \textnormal{$v$ is a multiplier of $A(G)$} \}
\]
of $B (A(G))$, where  $B (A(G))$ denotes the algebra of all bounded linear 
operators from $A(G)$ to $A(G)$. 
Throughout this paper we will generally use $v$ in place of the 
operator $M_{v}$ and we will
write $\parallel v\parallel _{MA(G)}$ to represent the norm of $M_{v}$ in  $B (A(G))$.

\begin{defn}

Let $\cQ_{M}(G)$ be the completion of $L^1(G)$ with respect to the norm 
\[\|f\|_{\cQ_{A_M(G)}(G)} =\sup \left\{\int_G f(x) v(x) dx ~|~ v\in MA(G), \|v\|_{MA(G)}\leq1\right\}.\]

Then $\cQ_{M}(G)^*=MA(G)$. 
\end{defn}

Let $G$ be a locally compact group and let $VN(G)$ denote its group 
von Neumann algebra. The duality
\[
  A(G) = VN(G)_*\]
equips $A(G)$ with a natural operator space structure. With this 
operator space structures we can define the 
$cb$\emph{-multiplier algebra} of $A(G)$ to be
\[
  {M}_{cb}A(G) := CB(A(G)) \cap {M}(A(G)),
\]
where $CB(A(G))$ denotes the algebra of all completely bounded 
linear maps from $A(G)$ into itself. 
We let $\parallel  v\parallel _{cb}$ denote the $cb$-norm of 
the operator $M_{v}$. 
It is well known that ${M}_{cb}(A(G))$ is a closed subalgebra 
of $CB(A(G))$ and 
is thus a (completely contractive) Banach algebra
with respect to the norm $\parallel \cdot \parallel _{cb}$. 

\begin{defn} 

Let $\cQ_{cb}(G)$ be the completion of $L^1(G)$ with respect to the norm 
\[\|f\|_{\cQ_{cb}(G)} =\sup \left\{\int_G f(x) v(x) dx ~|~ v\in M_{cb}A(G), \|v\|_{cb}\leq1\right\}.\]

Then $\cQ_{cb}(G)^*=M_{cb}A(G)$.

\end{defn}

It is known that in general, 
\[A(G)\subseteq B(G) \subseteq {M}_{cb}(A(G)) \subseteq {M}(A(G))\] 
and that for $v \in A(G)$
\[\parallel v\parallel _{A(G)}=\parallel v\parallel _{B(G)}\ge 
\parallel v\parallel _{cb}\ge \parallel v\parallel _{M}.\]

In case $G$ is an amenable group, we have 
\[B(G)={M}_{cb}(A(G)) = {M}(A(G)) \] 
and that 
\[\parallel v\parallel _{B(G)} = \parallel v\parallel _{cb} = \parallel v\parallel _{M}\]
for any $v\in B(G)$.

\begin{defn}  Given a locally compact group $G$ let 
\[A_{M}(G)\stackrel{def}{=}A(G)^{-\parallel \cdot \parallel _{M } }\subseteq M(A(G)).\] 
and
\[A_{cb}(G)\stackrel{def}{=}A(G)^{-\parallel \cdot \parallel _{cb } }\subseteq M_{cb}(A(G)).\] 

\end{defn}

\begin{rem} The algebra $A_{cb}(G)$ was introduced by the first author in \cite{For} where it was denoted by $A_0(G)$. In that paper we show that in the case of 
$\mathbb{F}_2$, the free group on two generators, $A_{cb}(G)$ shares many of the properties characteristic of the Fourier-algebra of an amenable group. In particular, the algebra 
$A_{cb}(\mathbb{F}_2)$ is known to have a bounded approximate identity. The  locally compact groups $G$ for which $A_{cb}(G)$ has a bounded approximate identity are called \textit{weakly amenable groups}. All amenable groups are weakly amenable, but many classical non-amenable groups such as $\mathbb{F}_2$ and $SL(2,\mathbb{R})$ are weakly amenable. We say that the locally compact group $G$ is M-weakly amenable  if there is an approximate identity 
$\{u_{\alpha }\}_{\alpha \in I}$ in $A(G)$ that is bounded in the norm $\parallel \cdot \parallel _{M}$.

\end{rem}

\begin{rem} Let $\cA(G)$ denote either $A_{M}(G)$ or $A_{cb}(G)$. Consider the following map and its adjoints:

\smallskip
\quad  $i_{\cA}:  A(G)\rightarrow \cA(G)$    
 
\smallskip
\quad  $i_{\cA}^*:  \cA(G)^*\rightarrow VN(G)$     

\smallskip
\quad  $i_{\cA}^{**}:  VN(G)^*\rightarrow \cA(G)^{**}$,  

\noindent
where $i_{\cA}$ denotes the inclusion map. Since $i_{\cA}$ has dense range, $i_{\cA}^*$ is injective and as such is invertible
with inverse ${i_{\cA}^*}^{-1}$ on $Range (i_{\cA}^*)$. It is easy to see that $i_{\cA}^{*}$ is simply the restriction map. That is 
\[i_{\cA}^{*}(T)=T_{|_{A(G)}}.\]

It will also be useful to view all of the above maps as embeddings. That is, when $G$ is  
nonamenable $\cA(G)^*$  can be viewed as a proper subset of $VN(G)$ and $VN(G)^{*}$ as a proper subset of $\cA(G)^{**}$, of course with different norms. 

\end{rem}

\begin{defn} We let
\[X_{M}^{\lambda}(G)=\{T\in \cQ_{M}(G)~|~\langle u,T \rangle  =0~\textnormal{for every~} u\in A_{M}(G)\}\]
and let 
\[B_{M}^{\lambda}(G)=\overline{A_{M}(G)}^{\textnormal{weak}^*}\subseteq MA(G). \]

Let 
\[\cQ_{M}^{\lambda}(G)=\cQ_{M}(G)/X^{\lambda}_{M}.\]

\vspace{.07in}

 We let
\[X_{cb}^{\lambda}(G)=\{T\in \cQ_{A_{cb}(G)}(G)~|~\langle u,T \rangle  =0~\textnormal{for every~} u\in A_{cb}(G)\}\]
and let 
\[B_{cb}^{\lambda}(G)=\overline{A_{cb}(G)}^{\textnormal{weak}^*}\subseteq M_{cb}A(G). \]

Let 
\[\cQ_{cb}^{\lambda}(G)=\cQ_{cb}(G)/X^{\lambda}_{cb}.\]

\end{defn}
\begin{rem} $\cQ_{M}(G)$ and $\cQ_{cb}(G)$ are natural analogs of the group $C^*$-algebra $C^*(G)$. Similarly, $\cQ_{M}^{\lambda}(G))$ and $\cQ_{cb}^{\lambda}(G))$ 
are natural anaolgs of the reduced group $C^*$-algebra $C^*_{\lambda}(G)$. In particular, we have that 
\[(\cQ_{M}^{\lambda}(G))^*=B_{M}^{\lambda}(G)\]
and 
\[(\cQ_{cb}^{\lambda}(G))^*=B_{cb}^{\lambda}(G).\]

Moreover, we can identify $\cQ_{M}^{\lambda}(G)$ and $\cQ_{A_{cb}(G)}^{\lambda}(G)$ with the closed  subspaces of $A_{M}(G)^*$ and $A_{cb}(G)^*$ respectively generated by 
\[\{\phi_f~|~f \in L^1(G)\}\]
where 
\[\langle \phi_f, u \rangle  =\int_G f(x)u(x) ~dx.\]
As such going forward, we will view $\cQ_{M}^{\lambda}(G)$ and $\cQ_{cb}^{\lambda}(G)$ as norm closed subspaces of $A_{M}(G)^*$ and $A_{cb}(G)^*$ respectively. Moreover, we will generally identify the functional $\phi_f$ with the function $f\in L^1(G)$ for notational convenience.

Finally, we let 

\[\cQ_{M}^{\delta}(G)\stackrel{def}{=}\overline{\textnormal{span}\{\phi_x \mid x \in G\}}^{-\|\cdot\|_{A_M(G)^*}}\subseteq A_M(G)^*\] 
and 
\[\cQ_{cb}^{\delta}(G)\stackrel{def}{=}\overline{\textnormal{span}\{\phi_x \mid x\in G\}}^{-\|\cdot\|_{A_{cb}(G)^*}}\subseteq A_{cb}(G)^*\] 

where $\phi_x(v)=v(x)$ for each $x \in G$ and each $v\in A_M(G)$ and $A_{cb}(G)$ resepctively. In particular, if $G$ is discrete, then $\cQ_{M}^{\lambda}(G)=\cQ_{M}^{\delta}(G)$ and 
$\cQ_{cb}^{\lambda}(G)=\cQ_{cb}^{\delta}(G)$. 

\end{rem}

\section{Invariant Subspaces of $A_{cb}(G)$ or $A_M(G)$ and their functorial properties.}

In this  section we will focus our attention on some important invariant subspaces of the algebras $A(G)$, $A_{cb}(G)$, $A_M(G)$ and their closed ideals.

Let $\cA(G)$ be any of the algebras $A(G)$, $A_{cb}(G)$ and $A_M(G)$. We note that while by definition 
\[UCB(\cA(G)) = \overline{\textnormal{span}\{v \cdot T\mid v\in \cA(G), T\in \cA(G)^* \}}^{-\|\cdot\|_{\cA^*}}\]
it turns out that we don't actually need the linear span. 

\begin{prop}\label{UCB} Let $\cA(G)$ be any of the algebras $A(G)$, $A_{cb}(G)$ and $A_M(G)$. Let $I$ be a closed ideal in $\cA(G)$ with $Z(I)$ a set of spectral synthesis. Then 
\[UCB(I) = \overline{\{v \cdot T\mid v\in \cA(G), T\in I^* \}}^{-\|\cdot\|_{I^*}}\]
    
\end{prop}

\begin{pf} Let 
\[j(Z(I))=\{u\in \cA(G) ~|~ supp(u) \textnormal{~is compact and disjoint from~} Z(I)\}.\]
The assumption that $Z(I)$ a set of spectral synthesis tells us that $j(Z(I))$
is dense in I.     

Let $u_1,u_2\in j(Z(I))$ and $T_1,T_2\in I^*$. Then if $K_1=supp(u_1)$ and $K_2=supp(u_2)$ we have that $K=K_1\cup K_2$ is a compact set disjoint from $Z(I)$. 
By the regularity, we can also find a $u\in A(G)$ with compact support disjoint from Z(I) such that $u(x)=1$ if $x\in K$. In particular, $u\in j(Z(I))$. Moreover 
\[u\cdot(u_1\cdot T_1+u_2\cdot T_2)=u u_1\cdot T_1+u u_2\cdot T_2=u_1\cdot T_1+ u_2\cdot T_2.\]
Hence $j(Z(I))\cdot I^*$ is a subspace of $UCB(I)$. 

Next, we assume that $u\in I$ and $T\in I^*$. Choose $\{u_n\}\subset j(Z(I))$
such that $\lim\limits_{n\to \infty} \|u_n-u\|_{\cA(G)}=0$. Since 
\[\|u_n \cdot T - u\cdot T\|_{I^*} \leq \|u_n-u\|_{\cA(G)} \| T\|_{I^*}\]
we have that $j(Z(I))\cdot I^*$ is also dense in $UCB(I)$. 

\end{pf}

\begin{rem} Given an ideal $I \in \cA(G)$ and a $T\in I^*$ we say that $x\in G \setminus Z(I)$ is a support point for $T$ if for every neighbourhood $V$ of $x$, there exists a $u\in I$ with $supp(u)\subseteq V$ and $\langle T,u \rangle \not = 0$. We let $supp(T)$ denote the collection of all support points of $T$. We note that Proposition \ref{UCB} shows that $UCB(I)$ is the closure in $I^*$ of those $T\in I^*$ with compact support. 

\end{rem}

\begin{defn} Let $\cA(G)$ denote any of the algebras $A(G)$, $A_{cb}(G)$ or $A_M(G)$. Given a closed ideal $I$ in $\cA(G)$ we define 
\[UCB_c(I)=\{T\in I^*~|~T ~\textnormal{has compact support}\}\]
or equivalently
\[UCB_c(I)= \{T\in I^*~|~T=u\cdot T_1, u\in I ~\textnormal{has compact support and }T_1 \in I^* \}.\]
    
\end{defn}

\begin{defn} Let $X$ be a closed invariant subspace of $\cA^*$ containing  a character $\phi\in \Delta(\cA)$. We say that $m\in X^*$ is a topological invariant mean at $\phi$ on $X$ if 
\[\|m\|_{X^*}=\langle m,\phi \rangle =1\] 
and 
\[\langle m, u\cdot T \rangle  = \langle \phi, u \rangle  \langle m,T \rangle \]
for all $u\in \cA$ and $T \in X$. We denote the space of all topological invariant means at $\phi$ on $X$ by $TIM_{\cA}(X,\phi)$.

\end{defn}

\begin{defn} Let $\cA(G)$ be one of the algebras  $A(G)$, $A_{cb}(G)$ or $A_M(G)$. Let $x\in G$ define the isometry $L_x:\cA \to \cA$ by 
\[L_x(u)(y) = u(xy),\]
for each $y \in G$. 

\end{defn}

The following two propositions  are respectively \cite[Proposition 4.3]{FST} and \cite[Proposition 4.5]{FST}.

\begin{prop}  Let $\cA(G)$ be one of the algebras  $A(G)$, $A_{cb}(G)$ or $A_M(G)$. Let $X\subseteq \cA(G)^*$ be a closed invariant subspace of $\cA(G)$. Let $x\in G$. 
\begin{itemize} 
\item[i)] If $Y= L^*_x(X)$, then $Y$ is a closed invariant subspace of  $\cA(G)^*$ and
\[u\cdot L^*_x(T)=L_x^*(L_x(u) \cdot T)\]
for every $u\in \cA$ and $T \in X$. 

\item[ii)] Let $x\in G$. Then $L^*_x(\phi_e)=\phi_x$ where $e$ denotes the identity of $G$. 
\item[iii)] Let $m\in TIM_{\cA(G)}(X,\phi_e)$. If   $x\in G$, Then $\phi_x\in L_x^*(X)$ and
 $L^{**}_{x^{-1}}(m) \in TIM_{\cA(G)}(L^*_x(X),\phi_x)$. 
\end{itemize}
\end{prop}

\begin{prop}\label{UCB1} Let $\cA(G)$ be either $A_{cb}(G)$ or $A_M(G)$. Let $\cA(G)\cdot VN(G) = \{ u\cdot T: u\in \cA(G), T\in VN(G) \}$. Then

\smallskip
\begin{enumerate}

\item[i)] $\cA(G)\cdot VN(G)\subseteq UCB(A(G))$

\item[ii]  $i^*(v\cdot T)=v\cdot i^*(T)$ for each $v\in \cA(G), T\in \cA(G)^{*}$.

\item[iii)] $i^{*}(UCB(\cA(G)))\subseteq UCB(A(G))$.

\item[iv)] If $\cA$ has a bounded approximate identity, then $\cA(G)\cdot VN(G) \cdot VN(G) = UCB(A(G))$ 

\item[v)] $u \cdot T \in i^{*}(\cA(G)^{*})$ for each $u\in A(G), T\in VN(G)$.

\end{enumerate}

 \end{prop} 

 We can now extend the result in Proposition \ref{UCB1} iii) to other important invariant subspaces.

 \begin{prop} Let $\cA(G)$ be either $A_{cb}(G)$ or $A_M(G)$. Then 
   \begin{itemize}
     \item[i)] $i^*(AP(\cA(G))\subseteq AP(A(G))$. 
      \item[ii)] $i^*(WAP(\cA(G))\subseteq WAP(A(G))$. 
      \item[iii)] $i^*(\cQ^{\delta}_{cb}(G))\subseteq C^*_{\delta}(G)$ and $i^*(\cQ^{\delta}_{M}(G))\subseteq C^*_{\delta}(G)$. 
      \item[iv)] $i^*(\cQ^{\lambda}_{cb}(G)) \subseteq C^*_{\lambda}(G)$ and $i^*(\cQ^{\lambda}_{M}(G)) \subseteq C^*_{\lambda}(G)$. 
 \end{itemize}  
 \end{prop}
 \begin{pf}
\begin{itemize}
     \item[i)] Let $T\in AP(\cA(G))$. Let $\{u_n\}$ be a sequence in $A(G)$ with 
     $\|u_n\|_{A(G)}\leq 1$ for each $n \in \mathbb{N}$. Then $
     \|u_n\|_{\cA(G)}\leq 1$ for each $n \in \mathbb{N}$. It follows that there exists a subsequence $\{u_{n_k}\}$ such that $\{u_{n_k}\cdot T\}$ converges to some $T_1 \in AP(\cA(G))$. Since 
     \[i^*(u_{n_k} \cdot T) = u_{n_k} \cdot i^*(T)\]
     and since $i^*$ is norm continuous we have that $\{u_{n_k} \cdot i^*(T)\}$ converges to $i^*(T_1)$. 
      \item[ii)] Let $T\in WAP(\cA(G))$. . Let $\{u_{\alpha}\}_{\alpha  \in \Omega}$ be a net in  in $A(G)$ with 
     $\|u_{\alpha}\|_{A(G)}\leq 1$ for each $\alpha  \in \Omega$. Then as before 
     $\|u_{\alpha}\|_{\cA(G)}\leq 1$ for each $\alpha  \in \Omega$. It follows that there exists a subnet  $\{u_{\alpha_{\beta}}\}$ such that 
     $\{u_{\alpha_{\beta}}\cdot T\}$ converges weakly to some $T_1 \in WAP(\cA(G)$. 

     Let $m\in A(G)^{**}$. Then $i^{**}(m)\in \cA(G)^{**}$. It follows that 
    \begin{eqnarray*}   
   \lim\limits_{\beta}\langle m, u_{\alpha_{\beta}} \cdot i^*(T) \rangle  &=& \lim\limits_{\beta}\langle m, i^*(u_{\alpha_{\beta}} \cdot T) \rangle \\
   &=& \lim\limits_{\beta}\langle i^{**}(m), u_{\alpha_{\beta}} \cdot T \rangle \\
   &=& \langle i^{**}(m), T_1 \rangle \\
   &=&\langle m,  i^*(T_1) \rangle . \\
    \end{eqnarray*}
This shows that  $\{u_{n_k}\cdot i^*(T)\}$ converges weakly to $i^*(T_1)$ and hence $i^*(T)\in WAP(A(G))$. 

\item[iii)] This follows immediately since, abusing notation, $i^*(\phi_x)=\phi_x$ for every $x\in G$. 
     
  \item[iv)]  Again, this follows immediately since, once more abusing notation, 
  $i^*(f)=f$ for every $f\in L^1(G)$. 
 \end{itemize}  
 \end{pf}

 It is worth asking whether the containments in the previous proposition can ever be equalities. This is clearly the case if $G$ is amenable since $A(G)=A_M(G)=A_{cb}(G)$. In \cite[Lemma 4.7]{ForMiao} the first author together with T. Miao showed that if $i^{*}(UCB(\cA(G)))= UCB(A(G))$, then $G$ must be amenable. In fact, the following is true: 

 \begin{thm} Let $\cA(G)$ be either $A_{cb}(G)$ or $A_M(G)$. Let $X$ be a closed invariant subspace of $\cA(G)^*$ which contains $L^1(G)$.  If $i^*(X)$ is closed in $VN(G)$, then $G$ is amenable. 
     
 \end{thm}

 \begin{pf} 
 Assume that $i^*(X)$ is closed in $VN(G)$. Since $i^*(L^1(G))=L^1(G)$, we have that $C_{\lambda}^*(G)\subseteq i^*(X)$. Moreover, since $i^*$ is continuous and injective, if $i^*(X)$ is closed, the restriction of $i^*$ to $X$ is an isomorphism. In particular, if we let $u\in \cA(G)$, then $u$ defines a continuous linear fucntional on $C_{\lambda}^*(G)$. Hence $\cA(G)\subseteq B_{\lambda}(G)$ and the norm 
 on $\cA(G)$ is equivalent to the norm from $B(G)$. As such $G$ is amenable \cite{Losert}.     
 \end{pf}

\begin{rem} Let $\cA(G)$ be either $A_{cb}(G)$ or $A_M(G)$. Let $H$ be a closed   subgroup of $G$. Then the restriction map $R: \cA(G)\to \cA(H)$ is a
 contraction homomorphism. In general, we don't know if it is surjective. In fact, if $G=SL(2,\mathbb{R})$ and $H$ is a closed subgroup of $G$ isomorphic to $\mathbb{F}_2$, 
  then the restriction map $R:A_M(G)\to A_M(H)$ is not surjective \cite[Theorem 4.3]{BranFor}. However, there are cases when the restriction map is known to be surjective.
 For example, we say that $G\in [SIN]_H$ if there is a neighborhood base at the identity consisting of sets that are invariant under the inner automorphisms by elements in $H$.
 If  $G\in [SIN]_H$, then $R:\cA(G)\to \cA(H)$ is surjective \cite[Corollary 3.4]{BranFor}.

If $H$ is amenable, then $\cA(H)=A(H)$ and the surjectivity of $R$ follows from Herz's Extension Theorem \cite{Herz}. In this case $A(H)$ is isometrically isomorphic to $\cA(G)/I_{\cA(G)}(H)$. Moreover, in this case $R^*$ is a $*$-homomorphism from $VN(H)$ onto \[VN_H(G)=\overline{\textnormal{span}\{\phi_x ~|~x\in H\} }^{-w*}\]
which is a von Neumann algebra in $\cA(G)^*$. (See \cite[Lemma 3.1]{KanLau}). Moreover, 
\[VN_H(G)=I_{\cA(G)}(H)^{\perp}.\]

In general, we let
\[\cA^*_H(G)=\overline{\textnormal{span}\{\phi_x ~|~x\in H\} }^{-w*}=I_{\cA(G)}(H)^{\perp}\subseteq \cA(G)^*.\]

\end{rem}

\begin{defn} Let $G$ be a locally compact group, $H$ a closed subgroup and $C \geq 1$ be a constant. We say that 
$\mathcal{A}(H)$ is \textit{C-extendable} in $\mathcal{A}(G)$ if for every $v \in \mathcal{A}(H)$ and every $\epsilon > 0$, there 
exists a $w \in \mathcal{A}(G)$ such that $v= w_{|_H}$ and 
\[\|w\|_{ \mathcal{A}(G)} \leq (C+\epsilon) \|v\|_{\mathcal{A}(H)} .\]
Moreover, we also assume that if $v\in A(G)$, then we can choose $w\in A(G)$. 

We let 

\[Ext_{\mathcal{A}(G)}=\{H |  \mathcal{A}(H) \textnormal{~is $C$-extendable for some $C\geq 1$}\}.\]

\end{defn}

\begin{rem} $\cA(G)$ be either $A_{cb}(G)$ or $A_M(G)$. Assume that $H\in Ext_{\mathcal{A}(G)}$. Then $R:\cA(G)\to \cA(H)$ is surjective with kernel $I_{\cA(G)}(H)$ and hence $\cA(H)$ is isomorphic with $\cA(G)/I_{\cA(G)}(H)$. In particular, if $\mathcal{A}(H)$ is $C$-extendable in $\mathcal{A}(G)$, then for each $T\in \cA(H)^*$
\[1/C \|T\|_{\cA(H)^*}\leq \|R^*(T)\|_{\cA(G)^*}\leq \|T\|_{\cA(H)^*}\] and $R^*$ is an isomorphism of $\cA(H)^*$ onto $\cA^*_H(G)$. In case, $C=1$, the natural identification of  $\cA(H)$  with $\cA(G)/I_{\cA(G)}(H)$ is an isometric algebra isomorhism.    
\end{rem}

The following proposition was established by Kaniuth and Lau \cite[Lemma 3.2]{KanLau} for $A(G)$ for all closed subgroups $H$. 

\begin{prop}\label{Func1} Let $\cA(G)$ either $A_{cb}(G)$ or $A_M(G)$. Let $H\in Ext_{\cA(G)}$.  Let $R:\cA(G)\to \cA(H)$ be the restriction map. Then 
\begin{itemize}
    \item[i)] $R^{*}(UCB_c(\cA(H)) = UCB_c(\cA(G))\cap \cA^*_H(G)$ 
    \item[ii)] $R^{*}(UCB(\cA(H))= UCB(\cA(G))\cap \cA^*_H(G)$
    \item[iii)] $R^{*}(AP(\cA(H))=  AP(\cA(G))\cap \cA^*_H(G)$ 
    \item[iv)]  $R^{*}(WAP(\cA(H))= WAP(\cA(G))\cap \cA^*_H(G)$ 
\end{itemize}
 
\end{prop}

\begin{pf} 
\begin{itemize}
\item[i)] If $T\in (UCB_c(\cA(H))$, then there exists a $u\in A(G) \cap C_{00}(G)$ and a 
$T_1\in \cA(G)^*$ such that $T=u\cdot T_1$. Since $u\in A(G)$, there exists a $v\in A(G) \cap C_{00}(G)$ such that $R(v)=u$.

\begin{eqnarray*} 
\langle R^*(T), w \rangle  &=& \langle R^*(u \cdot T_1), w \rangle \\
&=& \langle u \cdot T_1, R(w) \rangle \\
&=& \langle  T_1, uR(w) \rangle \\
&=& \langle  T_1, R(vw) \rangle \\
&=& \langle  R^*(T_1), vw \rangle \\
&=& \langle v\cdot R^*(T_1), w \rangle .\\
\end{eqnarray*}

Hence $R^*(T)=v\cdot R^*(T_1)\in UCB_c(\cA(G))$. But we also have $R^*(T)\in \cA^*_H(G)$. 

Next assume that $T \in UCB_c(\cA(G))\cap \cA^*_H(G)$. In particular, there exists $T_1\in \cA(H)^*$ such that $T=R^*(T_1)$. 

Let $K=supp(T)$. Since $T\in UCB_c(\cA(G))$, $K$ is compact. As such there exists 
$u\in A(G)\cap C_{00}(G)$ such that $u(x)=1$ for each $x$ in some open neighborhood of $K$. It follows
that $u \cdot T =T$. Note also that $v=R(u)\in A(H)\cap C_{00}(H)$. In particular 
$v\cdot T_1\in (UCB_c(\cA(H))$. We claim that $R^*(v\cdot T_1)=T$.

Let $w\in \cA(G)$. Then 
\begin{eqnarray*} 
\langle R^*(v\cdot T_1), w \rangle  &=& \langle v\cdot T_1, R(W) \rangle  \\
&=& \langle  T_1, v R(w) \rangle \\
&=& \langle  T_1, R(uw) \rangle \\
&=& \langle  R^*(T_1), uw)\\
&=& \langle  u\cdot T, w \rangle \\
&=& \langle T, w \rangle .\\
\end{eqnarray*}
\item[ii)] Since $R^*$ is continuous, and $UCB_c(\cA(H))$ is dense in $UCB(\cA(H))$, it follows from i) that $R^{*}(UCB(\cA(H))\subseteq UCB(\cA(G))\cap \cA^*_H(G)$. 

Since  $H\in Ext_{\mathcal{A}(G)}$, $R^*$ is an isomorphism and hence $R^{*}(UCB(\cA(H))$ has closed range. It now follows from i) that $R^{*}(UCB(\cA(H))= UCB(\cA(G))\cap \cA^*_H(G)$

\item[iii)] Let $T\in AP(\cA(H))$. Let $\|u\|_{\cA(G)}\leq 1$. Then $\|R(u)\|_{\cA(H)}\leq 1$. Moreover, if $v=R(u)$, then for each $w\in \cA(G)$, we have 
\begin{eqnarray*}
    \langle R^*(v\cdot T), w \rangle  &=& \langle v \cdot T, R(w) \rangle  \\
    &=& \langle T, R(uw) \rangle \\
    &=& \langle R^*(T), uw \rangle \\
    &=& \langle u\cdot R^*(T), w \rangle \\
\end{eqnarray*}

It follows that 
\[\{u\cdot R^*(T)~|~\|u\|_{\cA(G)}\leq 1\}\subseteq R^*(\{v\cdot T~|~  \|v\|_{\cA(H)}\leq 1\}. \]

Since $\{v\cdot T~|~  \|v\|_{\cA(H)}\leq 1\}$ is relatively compact, so is $R^*(\{v\cdot T~|~  \|v\|_{\cA(H)}\leq 1\}$ and hence so is $\{u\cdot R^*(T)~|~\|u\|_{\cA(G)}\leq 1\}$. That is 
\[R^{*}(AP(\cA(H))\subseteq AP(\cA(G))\cap \cA^*_H(G).\]

Next suppose that $T\in AP(\cA(G))\cap \cA^*_H(G)$. It follows that there exists a $T_1\in \cA(H)^*$ such that $T=R^*(T_1)$. Let $\{u_n\}\subset \cA(H)$ with $\|u_n\|_{\cA(H)}\leq 1.$ For each $n\in \mathbb{N}$, choose $v_n \in \cA(G)$, such that $\|v_n\|_{\cA(G)^*}\leq C+1$ and $v_{n_{|_H}}=u_n$. Then from our calculation above we see that 
\[v_n \cdot T = R^*(u_n\cdot T_1).\]
Since $T\in AP(\cA(G))\cap \cA^*_H(G)$, it follows that there exists a subsequence 
$\{v_{n_k} \cdot T\}$ of $\{v_{n} \cdot T\}$, that converges in norm to some $S\in \cA^*_H(G)$. As $R^*$ is an isomorphism onto $\cA^*_H(G)$, there exists an $S_1 \in \cA(H)^*$ so that $\{u_{n_k} \cdot T_1\}$ converges in norm to $S_1$. This shows that $T_1 \in AP(\cA(H))$ and hence that $R^{*}(AP(\cA(H))=  AP(\cA(G))\cap \cA^*_H(G)$. 

\item[iv)] That $R^{*}(WAP(\cA(H))\subseteq  WAP(\cA(G))\cap \cA^*_H(G)$ follows in a manner similar to the proof of iii) above. 

Suppose that $T\in WAP(\cA(G))\cap \cA^*_H(G)$. It follows as before that that there exists a $T_1\in \cA(H)^*$ such that $T=R^*(T_1)$. Let $\{u_{\alpha}\}_{\alpha \in \Omega}\subseteq \cA(H)$ with $\|u_n\|_{\cA(H)}\leq 1.$ For each $\alpha \in \Omega$, choose $v_{\alpha} \in \cA(G)$, such that $\|v_{\alpha}\|_{\cA(G)}\leq C+1$ and $v_{\alpha_{|_H}}=u_{\alpha}$. 

Since $T\in WAP(\cA(G))\cap \cA^*_H(G)$, there exists a subnet $\{v_{\alpha_\beta} \cdot T\}$ of $\{v_{\alpha} \cdot T\}$ that converges weakly to some $S\in \cA(G)^*$. However,  since $ \cA^*_H(G)$ is weakly closed, we have that $S\in \cA^*_H(G)$. In particular, $\{v_{\alpha_\beta} \cdot T\}$ converges to $S$ in the $\sigma(\cA^*_H(G), \cA^*_H(G)^*)$ topology. Moreover $S=R^*(S_1)$ for some $S_1\in \cA(H)^*$. Since $R^*:\cA(H)^* \to \cA^*_H(G)^*$ is an isomorhpism it follows that 
$\{u_{\alpha_\beta} \cdot T_1\}$ converges to $S_1$ weakly in $\cA(H)^*$. Hence $T_1\in WAP(\cA(H))$. 

\end{itemize} 
\end{pf}

\textbf{Note:} In the previous proposition, we do not know if, in general, the assumption that $H \in Ext_{\mathcal{A}(G)}$ to obtain the equalities in the previous proposition. 

\begin{rem} Let $\cM(G)$ denote either $MA(G)$ or $M_{cb}A(G)$, and let $\cA(G)$ denote respectively $A_M(G)$ or $A_{cb}(G)$. Let $H$ be an open subgroup of $G$. Let $\cM(G)$ denote either For each $u\in \cM(H)$ 
we define $u^{\circ}$ by 

\begin{equation*}
u^{\circ} = \left\{
\begin{array}{rl}
u(x) & \text{if } x \in H,\\
0 & \text{if } x \not \in H.\\
\end{array} \right.
\end{equation*}

Then $u^{\circ}\in \cM(G)$   with $\|u^{\circ}\|_{\cM(G)}=\|u\|_{\cM(H)} $.   

Define  $\Gamma:\cM(H)\to \cM(G)$  by $\Gamma(u)=u^{\circ}$.  
Then $\Gamma $ is an isometric algebra isomorphism from  $\cM(H)$  onto the ideal  
$1_H \cM(G)$  of $\cM(G)$. Moreover, if $u\in \cA(H)$, 
 then $\Gamma(u)=u^{\circ}\in \cA(G)$  

Note that if $T\in \cA(G)^*$ and $w \in \cA(G)$, then 
\begin{eqnarray*}
    \langle R^*\circ \Gamma^*(T), w \rangle &=& \langle  \Gamma^*(T), R(w) \rangle \\
    &=& \langle  T, \Gamma(R(w)) \rangle \\
    &=& \langle T, 1_H w \rangle \\
    &=& \langle 1_H \cdot T, w \rangle .\\
\end{eqnarray*}
That is $R^*\circ \Gamma^*(T)= 1_H \cdot T$. 

If $T\in \cA(G)^*$, $v\in \cA(G)$ and $u=R(v)$, then for $w\in \cA(H)$, we have 
\begin{eqnarray*} 
\langle \Gamma^*(v\cdot T), w \rangle &=& \langle v\cdot T, w^{\circ} \rangle \\
&=& \langle  T, v w^{\circ} \rangle \\
&=& \langle  T, (uw)^{\circ} \rangle \\
&=& \langle  \Gamma^*(T), uw \rangle \\
&=& \langle  u\cdot\Gamma^*(T), w \rangle \\
\end{eqnarray*}
Hence $\Gamma^*(v\cdot T)=R(v)\cdot \Gamma^*(T)$. 

Finally, with $H$ open, we can define a projection $P_H$ from $\cA(G)$ onto the ideal 
$I_{\cA(G)}(G\setminus H)$ by $P_H(u)=1_H u$.The kernel of $P_H$ is the ideal $I_{\cA(G)}(H)$. Moreover, $P_H=\Gamma \circ R$. 

\end{rem}

\begin{prop} Let $\cA(G)$ denote either  $A_M(G)$ or $A_{cb}(G)$. Let $H$ be an open subgroup of $G$. Then 
\begin{itemize}
    \item[i)] $\Gamma^*(\cA(G)^*)=\cA(H)^*$ and $\Gamma^*_{|_{1_H \cdot \cA(G)^*}}$ is an isometric isomorphism from $1_H \cdot \cA(G)^*$ onto $\cA(H)^*$
    \item[ii)] $\Gamma^*(UCB(\cA(G)))=UCB(\cA(H))$ and $\Gamma^*_{|_{1_H \cdot UCB(\cA(G))}}$ is an isometric isomorphism 
    from $1_H \cdot UCB(\cA(G))$ onto $UCB(\cA(H))$. 
    
     \item[iii)] $\Gamma^*(AP(\cA(G)))=AP(\cA(H))$ and $\Gamma^*_{|_{1_H \cdot AP(\cA(G))}}$ is an isometric isomorphism from $1_H \cdot AP(\cA(G))$ onto $AP(\cA(H))$.
     \item[iv)]  $\Gamma^*(WAP(\cA(G)))=WAP(\cA(H))$ and $\Gamma^*_{|_{1_H \cdot WAP(\cA(G))}}$ is an isometric isomorphism from $1_H \cdot WAP(\cA(G))$ onto $WAP(\cA(H))$.
     \item[V)]  $\Gamma^*(Q^{\delta}_M(G))=Q^{\delta}_M(H)$ and 
     $\Gamma_M^*{_{|_{1_H \cdot Q^{\delta}_M(G)}}}$ is an isometric 
     isomorphism from $1_H\cdot Q^{\delta}_M(G)$ onto $Q^{\delta}_M(H)$. 
     \item[vi)]  $\Gamma^*(Q^{\delta}_{cb}(G))=Q^{\delta}_{cb}(H)$ and $\Gamma^*_{|_{1_H\cdot Q^{\delta}_{cb}(G)}}$ is an isometric isomorphism from $1_H\cdot Q^{\delta}_{cb}(G)$ onto $Q^{\delta}_{cb}(H)$.
     \item[vii)]  $\Gamma^*(Q^{\lambda}_M(G))=Q^{\lambda}_M(H)$ and $\Gamma^*_{|_{1_H\cdot Q^{\lambda}_M(G)}}$ is an isometric isomorphism from $1_H\cdot Q^{\lambda}_M(G)$ onto $Q^{\lambda}_M(H)$.
     \item[viii)]  $\Gamma^*(Q^{\lambda}_{cb}(G))=Q^{\lambda}_{cb}(H)$ and $\Gamma^*_{|_{1_H\cdot Q^{\lambda}_{cb}(G)}}$ is an isometric isomorphism from $1_H\cdot Q^{\lambda}_{cb}(G)$ onto $Q^{\lambda}_{cb}(H)$.
\end{itemize}
    
\end{prop}

\begin{pf} 
\begin{itemize} 
\item[i)] We first show that $\Gamma^*$ is surjective. Assume that $T\in \cA(H)^*$. Let $S=R^*(T)$. We claim that $T=\Gamma^*(S)$. Let $u\in \cA(G)$. Then 
\begin{eqnarray*} 
\langle \Gamma^*(S), u \rangle &=& \langle S, \Gamma(u) \rangle \\
&=& \langle R^*(T), u^{\circ} \rangle \\
&=& \langle T, R(u^{\circ}) \rangle \\
&=& \langle T, u \rangle \\
\end{eqnarray*} 

Next we show that $ker(\Gamma^*)=1_{G\setminus H}\cdot \cA(G)^*$. Let 
$1_{G\setminus H}\cdot T\in 1_{G\setminus H} \cdot \cA(G)^*$ and let $u\in \cA(H)$. Then 
\begin{eqnarray*} 
\langle \Gamma^*(1_{G\setminus H} \cdot T), u \rangle &=&\langle 1_{G\setminus H} \cdot T, \Gamma(u) \rangle \\
&=&\langle 1_{G\setminus H} \cdot T, u^{\circ} \rangle \\
&=& \langle T, 1_{G\setminus H} u^{\circ} \rangle \\
&=& 0.\\
\end{eqnarray*}

Assume that $T\in \cA(G)^*$, but $T\not \in 1_{G\setminus H} \cdot \cA(G)^*$. Then
there exists $x\in H \cap supp(T)$. In particular, since $H$ is an open neighbourhood 
of $x$, there exists $v\in \cA(G)$ such that $v(x)\not =0$ with $supp(v)\subseteq H$ 
and $\langle T, v \rangle \not = 0$. Let $u= R(v)$. Then $u^{\circ}=v$ so 
\[\langle \Gamma^*(T),u \rangle =\langle T,\Gamma(u) \rangle =\langle T,u^{\circ} \rangle =\langle T,v \rangle \not = 0. \]

Finally, we show that $\Gamma^*_{|_{1_H \cdot \cA(G)^*}}$ is an isometry. Clearly, $\|\Gamma^*\|_{B(\cA(G)^*,\cA(H)^*)}\leq 1$. 

Let $T \in 1_H \cdot \cA(G)^*$  such that $\|T\|_{\cA(G)^*}=1$. Let $\epsilon > 0$. Then  there exists $v\in \cA(G)$ with $\|v\|_{\cA(G)}=1 $ such that 
\[|\langle T,v \rangle | > 1-\epsilon.\] Note that since  $T \in 1_H \cdot \cA(G)^*$, we have that $T=1_H \cdot T$ so that if $w=1_H v$, then $\|w\|_{\cA(G)}\leq 1 $ and 
\[|\langle T,v \rangle |=|\langle T,w \rangle | > 1-\epsilon.\]
Let $u=R(v)$. Then $\|u\|_{\cA(H)}\leq 1 $ and 
\[|\langle \Gamma^*(T),u \rangle |=|\langle T,\Gamma(u) \rangle |=|\langle T,w \rangle |  > 1-\epsilon.\]
\item[ii)] We need only show that $\Gamma^*(UCB(\cA(G)))=UCB(\cA(H))$. 

The first observation we will make is that for any $T\in \cA(G)^*$, we can write 
\[T=1_H\cdot T + 1_{G\setminus H}\cdot T\]
and that $\Gamma^*(T)=\Gamma^*(1_H\cdot T)$. Moreover 
\[UCB(\cA(G))=1_H\cdot UCB(\cA(G)) \oplus 1_{G\setminus H}\cdot UCB(\cA(G)) \]

Assume that $T=v \cdot T_1 \in UCB(\cA(G))$ with $v\in \cA(G)$. Let $u=R(v)$. 
Then $\Gamma^*(T)=u\cdot \Gamma^*(T_1)\in UCB(\cA(H))$. Since $\Gamma^*$ is 
continuous and $\{u\cdot T~|~ T \in \cA(G)^*, u \in \cA(G)\}$ is dense in
$UCB(\cA(G))$, we get that $\Gamma^*(UCB(\cA(G)))\subseteq UCB(\cA(H))$. 

Next assume that $T= u\cdot T_1$ where $T_1 \in \cA(H)^*, u \in \cA(H)$. Then there exists an $S\in 1_H \cdot\cA(G)^*$ such that $T_1=\Gamma^*(S)$. If $v=u^{\circ}$, then 
\[\Gamma^*(v\cdot S) = R(v)\cdot \Gamma^*(S)=u\cdot T_1. \] 
Hence, $\Gamma^*: 1_H \cdot UCB(\cA(G))\to UCB(\cA(H)$ has dense range. But as $\Gamma^*$ is isometric on $1_H \cdot UCB(\cA(G))$, we have that $\Gamma^*(UCB(\cA(G)))=UCB(\cA(H))$.

\item[iii)] It is easy to see that $T\in AP(\cA(G))$ if and only if both $1_H\cdot T$ and $1_{G\setminus H}\cdot T$ are in $AP(\cA(G)$. Let $T\in 1_H \cdot AP(\cA(G))$. Let $S=\Gamma^*(T)$. Now 
\[\{u\cdot S|\|u\|_{\cA(H)}\leq 1\} = \Gamma^*(\{u^{\circ}\cdot T|\|u\|_{\cA(H)}\leq 1\}).\]
But $\{u^{\circ}\cdot T|\|u\|_{\cA(H)}\leq 1\})$ is relatively compact. Hence,  
$\{u\cdot S|\|u\|_{\cA(H)}\leq 1\}$ is also relatively compact and $S\in AP(\cA(H))$.

For the converse, suppose that $S\in AP(\cA(H))$ and choose $T\in 1_H\cdot \cA(G)^*$, so that $S=\Gamma^*(T)$. Since $1_H\cdot \cA(G)^*$ is invariant, $\Gamma^*$ maps 
$\{v\cdot T|\|v\|_{\cA(G)}\leq 1\}$ isometrically onto $\{R(v)\cdot S|\|v\|_{\cA(G)}\leq 1\}$. Since  $\{R(v)\cdot S|\|v\|_{\cA(G)}\leq 1\}$ is relatively compact, so is $\{v\cdot T|\|v\|_{\cA(G)}\leq 1\}$.

\item[iv)] This follows in a similar manner to iii) above. 

\item[v)] This follows from i) since, abusing notation, 
\begin{equation*}
\Gamma^*(\phi_x) = \left\{
\begin{array}{rl}
\phi_x & \text{if } x \in H,\\
0 & \text{if } x \not \in H.\\
\end{array} \right.
\end{equation*}
for every $x\in G$. 
\item[vi)] See v). 
\item[vii)] This follows from i) since, abusing notation, $\Gamma^*(\phi_f)= \phi_{f_{|_H}}$ and that $L^1(H)$ is isometrically isomorphic with $1_H L^1(G)$.
\item[viii)] See vii). 
\end{itemize} 
\end{pf}

\section{Containment Results} 

In this section, we will establish various containment results for invariant subspaces of $A(G)$, $A_{cb}(G)$, and $A_M(G)$ and their ideals.

\begin{lem}\label{UCB2} Let $\cA(G)$ denote any of the algebras $A(G)$, $A_{cb}(G)$ or $A_M(G)$. Let $I$ be a non-zero closed ideal in $\cA(G)$ with $E=Z(I)$. Then   
\[\overline{\textnormal{span}\{\phi_x \mid x\in G\setminus E\}}^{-\|\cdot\|_{I^*}}\subseteq UCB(I).\] 

If $G$ is discrete and $Z(I)$ is a set of spectral synthesis, then 
\[\overline{\textnormal{span}\{\phi_x \mid x\in G\setminus E\}}^{-\|\cdot\|_{I^*}}= UCB(I).\] 

\end{lem}

\begin{pf} Given any $x\in G\setminus E$, we can find an open neighbourhood $U$ of $x$ for which $U\cap E=\emptyset$. Then we can choose  $u\in A(G)\cap C_{00}(G)\subseteq \cA(G)$ 
so that $u(x)=1$ and $u(y)=0$ for every $y\not \in U$ and hence $u\in I$. If $v\in I$, then 
\[\langle u \cdot \phi_x, v \rangle  = \langle  \phi_x, vu \rangle =v(x)u(x)=v(x)=\langle \phi_x,v \rangle .\]
Hence $\phi_x= u \cdot \phi_x \in  UCB(I)$. It follows that $\overline{\textnormal{span}\{\phi_x \mid x\in G\setminus E\}}^{-\|\cdot\|_{I^*}}\subseteq UCB(I)$. 

Next assume that $G$ is discrete and that $E=Z(I)$ is a set of spectral synthesis. Let $T= u\cdot T_1\in I^*$ where $u\in I$ and $T_1\in I^*$. 

Let $x \in G\setminus E$. If 
\[J=\{v\in I\mid v(x)=0\},\]
then $J$ is a closed ideal in $I$ of co-dimension 1. Hence $J^{\perp}$ is one dimensional. Clearly $\phi_x\in J^{\perp}$. It follows that 
\[ J^{\perp} =\{ \alpha \phi_x \mid \alpha \in \mathbb{C}\}.\]
Next we note that if  $w=1_{\{x\}}$, then  $w\cdot T_1\in J^{\perp}$. It follows that $w\cdot T_1=\alpha \phi_x $ for some $\alpha \in \mathbb{C}$. 

Let $v\in c_{00}(G)\cap I$ have support $\{x_1,x_2,\ldots, x_k\}\subseteq G\setminus E$. Then $v=\sum\limits_{i=1}^k v(x_i) 1_{x_i}$ so that 
\[v\cdot T_1 = \sum\limits_{i=1}^k v(x_i) 1_{x_i} \cdot T_1= \sum\limits_{i=1}^k v(x_i) \alpha_i \phi_{x_i}\]

for some finite collection of scalars  $\{\alpha_1,\alpha_2,\ldots, \alpha_k\}$.

Finally since $E$ is a set of spectral synthesis, there exists a sequence 
$\{u_n\}\subset I$ with $u_n \in c_{00}(G)$ such that 
$\lim\limits_{n\to \infty}\| u_n-u\|_{\cA(G)} =0$. If $w \in I$ with $\|w\|_{\cA(G)}\leq 1$, then 
\begin{eqnarray*}
\mid \langle u\cdot T_1, w \rangle - \langle u_n \cdot T_1, w \rangle \mid&=& \mid \langle  T_1,wu \rangle - \langle  T_1,wu_n \rangle \mid\\
&=&\mid  \langle  T_1,w(u-u_n) \rangle \mid\\
&\leq& \|T_1\|_{I^*} \|u-u_n\|_{\cA(G).}\\
\end{eqnarray*} 

It follows that $\lim\limits_{n\to \infty}\| u_n\cdot T_1-u\cdot T_1\|_{I^*} =0$. However, 
$u_n\cdot T_1 \in \textnormal{span}\{\phi_x \mid x\in G\setminus E\}$ and hence we have that  $T=u\cdot T_1\in \overline{\textnormal{span}\{\phi_x \mid x\in G\setminus E\}}^{-\|\cdot\|_{I^*}}$.

 We get that 
\[ UCB(I) \subseteq \overline{\textnormal{span}\{\phi_x \mid x\in G\setminus E\}}^{-\|\cdot\|_{I^*}}.\]

\end{pf} 

Since the empty set is a set of Spectral synthesis for all of $A(G)$, $A_{cb}(G)$ and $A_M(G)$  the next result follows immediately. 

\begin{cor}\label{discrete1} Let $G$ be a locally compact group. Then $\cQ_{cb}^{\delta}(G)\subseteq UCB(A_{cb}(G))$ and $\cQ_{M}^{\delta}(G)\subseteq UCB(A_{M}(G))$. Moreover, if 
$G$ is discrete, then  $UCB(A_{cb}(G))=\cQ_{cb}^{\delta}(G)=\cQ_{cb}^{\lambda}(G)$ and 
$UCB(A_M(G))=\cQ_{M}^{\delta}(G)=\cQ_{M}^{\lambda}(G)$. 

\end{cor}  

\begin{rem} If $G$ is discrete, then  $UCB(A_{cb}(G))=\cQ_{cb}^{\lambda}(G)=\cQ_{cb}^{\delta}(G)$ and 
$UCB(A_{M}(G))=\cQ_{M}^{\lambda}(G)=\cQ_{M}^{\delta}(G)$.  It follows that 
\[UCB(A_{cb}(G))^*=B_{cb}^{\lambda}(G)\]
and 
\[UCB(A_{M}(G))^*=B_{M}^{\lambda}(G)\]
respectively. Moreover, the normal algebra structure of both $B_{cb}^{\lambda}(G)$ and $B_{M}^{\lambda}(G)$ agrees with the multiplication on $UCB(A_{cb}(G))^*$ and $UCB(A_{M}(G))^*$ respectively that is inherited from being the dual of a topologically introverted subspace. 

If $I$ is a closed ideal in $A_{cb}(G)$, and if $G$ is discrete and $E=Z(I)$ is a set of spectral synthesis, then we have seen that 
\[\overline{\textnormal{span}\{\phi_x \mid x\in G\setminus E\}}^{-\|\cdot\|_{I^*}}= UCB(I).\] 
In this case, 
\[UCB(I)^*=\overline{I}^{-w^*}\subseteq B_{cb}^{\lambda}(G).\]
Similarly, if $I\subset A_M(G)$, then 
\[UCB(I)^*=\overline{I}^{-w^*}\subseteq B_{M}^{\lambda}(G)\]
and if $I\subset A(G)$, then 
\[UCB(I)^*=\overline{I}^{-w^*}\subseteq B^{\lambda}(G).\]

\end{rem}

\begin{prop} Let $\cA(G)$ denote any of the algebras $A(G)$, $A_{cb}(G)$ or $A_M(G)$. Let $I$ be a non-zero closed ideal in $\cA$. If $I$ has a bounded approximate identity, then $WAP(I)\subseteq UCB(I)$. 

\end{prop} 

\begin{pf} Assume that I has a bounded approximate identity  $\{u_{\alpha}\}_{\alpha\in \Omega}$. Let $T\in WAP(I)$. Then for $u \in I$, 
\[\lim\limits_{\alpha\in \Omega} \langle  u_{\alpha} \cdot T, u \rangle  =  \lim\limits_{\alpha\in \Omega}\langle   T, u_{\alpha}u \rangle =\langle T, u \rangle .\]
That is $\{u_{\alpha}\cdot T\}_{\alpha\in \Omega}$ converges in the weak-$^*$ topology on $I^*$ to $T$. But since $T\in WAP(I)$ and since 
$\{u_{\alpha}\}_{\alpha\in \Omega}$ is bounded, there exists a subnet $\{u_{\alpha_{\beta}}\}_{\beta\in \Omega_1}$ such that 
$\{u_{\alpha_{\beta}}\cdot T \}_{\beta\in \Omega_1}$ converges weakly in $I^*$ to some $T_0 \in I^*$. But this means that $T=T_0$ so $T$ is in the weak closure of $\{u\cdot T| u\in I\}$ and hence in the weak closure of $UCB(I)$. But $UCB(I)$ is weakly closed so $T\in UCB(I)$. 
\end{pf}

\begin{rem}\label{baicont} It follows immediately from the definitions that for any commutative Banach algebra $\cA$, we have that $AP(\cA)\subseteq WAP(\cA)$.  Hence if $\cA(G)$ denotes any of the algebras $A(G)$, $A_{cb}(G)$ or $A_M(G)$ and if $I$ is a non-zero closed ideal in $\cA(G)$ with a bounded approximate identity, then $AP(I)\subseteq UCB(I)$. 
\end{rem}

\begin{thm}\label{AP1} Let $\cA(G)$ denote any of the algebras $A(G)$, $A_{cb}(G)$ or $A_M(G)$. 
\begin{itemize}
    \item[i)] Let $I$ be a non-zero closed ideal in $\cA(G)$ with $UCB(I) \subseteq WAP(I)$, then $G$ is discrete. 
    \item[ii)] Let $I$ be a non-zero closed ideal in $\cA$ with $x \not \in Z(I)$. Then $\phi_{x_{|_I}} \in AP(I)$. In particular, if $G$ is discrete and $Z(I)$ is a set of spectral synthesis, then 
$UCB(I)\subseteq AP(I)$. 
     \item[iii)] If $G$ is discrete, $Z(I)$ is a set of spectral synthesis, and $I$ has a bounded approximate identity, then $UCB(I)=AP(I)$.
\end{itemize}

\end{thm} 

\begin{pf} 
\begin{itemize} 
\item[i)] This is \cite[Corollary 4.19]{FST}.
\item[ii)] 
Let $\{u_{\alpha}\}_{\alpha\in \Omega}$ be a net in $I$  such that $\|u_{\alpha}\|_{\cA(G)}\leq 1$, we find a subnet $\{u_{\alpha_{\beta}}\}_{\beta\in \Omega_1}$ such  that $\lim\limits_{\beta\in \Omega_1}u_{\alpha_{\beta}}(x) =c\in \mathbb{C}$. It follows that if $u \in I$ with $\|u\|_{\cA(G)}\leq 1$, then 
\begin{eqnarray*} 
\mid \langle c \phi_{x_{|_I}}, u \rangle  - \langle u_{\alpha_{\beta}} \cdot  \phi_{x_{|_I}}, u \rangle \mid &=& \mid \langle \phi_{x_{|_I}}, cu \rangle  - \langle \phi_{x_{|_I}}, u_{\alpha_{\beta}}  u \rangle \mid \\
&=& \mid c u(x) - u_{\alpha_{\beta}}(x) u(x)\mid \\
&\leq & |c-u_{\alpha_{\beta}}(x)\mid \\
\end{eqnarray*}
This shows that $\lim\limits_{\beta\in \Omega_1} \|u_{\alpha_{\beta}} \cdot  \phi_{x_{|_I}}- c  \phi_{x_{|_I}}\|_{I^*}=0$. Hence 
$\phi_{x_{|_I}} \in AP(I)$. 

The fact that $UCB(I)\subseteq AP(I)$ follows from Lemma \ref{UCB2}

\item[iii)] This follows from ii) and Remark \ref{baicont}.

\end{itemize} 
\end{pf} 

\begin{cor} Let $\cA(G)$ denote either of the algebras  $A_{cb}(G)$ or $A_M(G)$.  Then $Q_{M}^{\delta}(G) \subset AP(\cA(G))$ and $Q_{cb}^{\delta}(G) \subset AP(\cA(G))$ . In particular,  $G$ is discrete if and only if $UCB(\cA(G))\subseteq AP(\cA(G))$. Moreover, if $\cA(G)$  has a bounded approximate identity, then $G$ is discrete if and only if $UCB(\cA(G))=AP(\cA(G))$.
\end{cor}

\begin{prop} Let $G$ be a locally compact group. Then 
\begin{itemize}
\item[i)] $\cQ^{\lambda}_{cb}(G)\subseteq UCB(A_{cb}(G))$
\item[ii)] $\cQ^{\lambda}_{cb}(G)\subseteq WAP(A_{cb}(G))$
\item[iii)] $\cQ^{\lambda}_{cb}(G)= UCB(A_{cb}(G))$ if and only if $G$ is discrete.
\item[iv)] $\cQ^{\lambda}_{M}(G)\subseteq UCB(A_{M}G)$. 
\item[v)] $\cQ^{\lambda}_{cb}(G)\subseteq WAP(A_{M}(G))$.
\item[vi)] $\cQ^{\lambda}_{M}(G)= UCB(A_{M}G)$ if and only if $G$ is discrete.
\end{itemize}
    
\end{prop}
\begin{pf}
 \begin{itemize}
 \item[i)] Since $C_{00}(G)$ is dense in $L^1(G)$ and hence in $\cQ^{\lambda}_{cb}(G)$, it is enough to show that every $f\in C_{00}(G)$ is in $UCB(A_{cb}(G))$. 
 
 Assume that  $f\in C_{00}(G)$ with $supp(f)=K$. Then there exits a $u \in A_{cb}(G)$ such that $u(x)=1$ for each $x\in K$. In particular $u\cdot f= uf= f$. Hence $f\in UCB(A_{cb}(G))$. 
 \item[ii)] Let $\{u_{\alpha}\}_{\alpha\in \Omega}$ be a net in $A_{cb}(G)$ with $\|u_{\alpha}\|_{A_{cb}(G)}\leq 1$ for each $\alpha\in \Omega$. Since 
 \[\{u_{\alpha}\}_{\alpha\in \Omega}\subseteq B^{\lambda}_{cb}(G)=\cQ^{\lambda}_{cb}(G)^*,\]
by passing to a subnet if necessary we can assume that  $\{u_{\alpha}\}_{\alpha\in \Omega}$ converges in the $\sigma(B^{\lambda}_{cb}(G),\cQ^{\lambda}_{cb}(G)) $ topology to some $u \in B^{\lambda}_{cb}(G)$.

Now consider $m\in A_{cb}(G)^{**}$. Let $v=m_{|_{\cQ^{\lambda}_{cb}(G)}}$. If $f\in \cQ^{\lambda}_{cb}(G)$, then 
\begin{eqnarray*}
\lim\limits_{\alpha\in \Omega}\langle m, u_{\alpha}\cdot f \rangle &=&\lim\limits_{\alpha\in \Omega}\langle v, u_{\alpha}\cdot f  \rangle \\
&=& \lim\limits_{\alpha\in \Omega} \int_G v(x) u_{\alpha}(x) f(x) dx\\
&=& \lim\limits_{\alpha\in \Omega}\langle u_{\alpha}, v\cdot f  \rangle  \\
&=& \langle u,v\cdot f \rangle \\
&=& \langle v,u \cdot f \rangle \\
&=&\langle m, u\cdot f \rangle .\\
\end{eqnarray*} 

That is $\{u_{\alpha}\cdot f\}_{\alpha\in \Omega}$ converges weakly to $u\cdot f$. This shows that $f\in WAP(A_{cb}(G))$. 
\item[iii)] Assume that $\cQ^{\lambda}_{cb}(G)= UCB(A_{cb}(G))$. Then by ii) above we have that $UCB(A_{cb}(G))\subseteq WAP(A_{cb}(G))$. In particular, Theorem \ref{AP1} i) shows that $G$ is discrete.

The converse is Corollary \ref{discrete1}. 

The proofs of iv), v) and vi) are similar to those of the corresponding statements about $A_{cb}(G)$. 
\end{itemize} 
\end{pf}

\begin{rem} If $\cA(G)$ is any of $A(G)$, $A_M(G)$ or $A_{cb}(G)$, then the question of when $\cA(G)=WAP(\cA(G))$ is equivalent to the question of identifying groups $G$ such that $\cA(G)$ is Arens regular. For $G$ amenable, it is known that $G$ must be finite. In general, while it is strongly suspected that  $\cA(G)=WAP(\cA(G))$ if and only if $G$ is finite,  surprisingly, even for $A(G)$, this is still not known to be true. However, we have recently shown that if  $\cA(G)=WAP(\cA(G))$, then $G$ must be discrete \cite[Theorem 5.1]{FST}. 

Granirer has shown that  $UCB(A(G))=VN(G)$ if and only if $G$ is compact. This motivates us asking: If $\cA(G)$ is either $A_M(G)$ or $A_{cb}(G)$, is $UCB(\cA(G))=\cA(G)^*$ if and only if $G$ is compact. While we strongly expect the answer to this question to be yes, we are unable to show this in general. However, we can show that if $UCB(\cA(G))=\cA(G)^*$, then $G$ has a compact open subgroup, and every amenable subgroup is finite. 

\end{rem} 

\begin{thm} Let  $\cA(G)$ be either $A_M(G)$ or $A_{cb}(G)$. If $UCB(\cA(G))=\cA(G)^*$, then $G$ has a compact open subgroup, and every closed amenable subgroup is finite. In particular, if $G$ is almost connected, then $G$ is compact.  

\end{thm}

\begin{pf} Assume that $H$ is a closed amenable subgroup. Then $H\in Ext_{\cA(G)}$. 
It follows from Proposition \ref{Func1} that 
$R^{*}(UCB(A(H))= UCB(\cA(G))\cap VN_H(G)$. However, if we assume that 
$UCB(\cA(G))=\cA(G)^*$, we have that $R^{*}(UCB(A(H))=VN_H(G)$. This in turn means that $UCB(A(H))=VN(H)$. It then follows from Granirer's result that $H$ is compact. 

To see that $G$ has an open subgroup we note that if the connected component $G_0$ is not amenable, then $G_0$ contains the free group $\mathbb{F}_2$ on two generators. In particular, $G_0$ contains an infinite non-compact abelian group which is impossible from our argument above. As a consequence it must be the case that $G_0$ is compact. From here we note that $G$ also has an open almost connected subgroup $G_1$. As $G_1/G_0$ is compact, $G_1$ must also be compact.  
\end{pf} 
\vskip 0.4 true cm

\begin{center}{\textbf{Acknowledgments}}
\end{center}
The authors wish to congratulate Professor Anthony To-Ming Lau on his forth coming 80th birthday and for his many and varied accomplishments over his distinguished career. In addition, the first author wishes to express his deep gratitude to Professor Lau for his exceptional mentorship, his unwavering support, and for the exemplary role model he has provided for all of us who have had the great pleasure to work with him.  \\ \\
\vskip 0.4 true cm

\vfill
\begin{tabbing}
{\it Second author's adddress\/}: \= Department of Mathematical and Statistical Sciences \kill 
{\it First author's adddress\/}:  \> Department of Pure Mathematics \\
                                  \> University of Waterloo \\
                                  \> Waterloo, Ontario \\
                                  \>  Canada N2L 3G1 \\[\medskipamount]
{\it E-mail\/}:                   \> {\tt beforres@uwaterloo.ca}\\[\bigskipamount]
{\it Second author's adddress\/}: \> Department of Pure Mathematics \\
                                  \> University of Waterloo \\
                                  \> Waterloo, Ontario \\
                                  \> Canada N2L 3G1 \\[\medskipamount]
{\it E-mail\/}:                   \> {\tt jmsawatz@uwaterloo.ca}\\[\medskipamount]
{\it \/}                     \>   \\[\bigskipamount]
{\it Third author's adddress\/}: \> Department of Pure Mathematics \\
                                  \> University of Waterloo \\
                                  \> Waterloo, Ontario \\
                                  \> Canada N2L 3G1 \\[\medskipamount]
{\it E-mail\/}:                   \> {\tt athamizhazhagan@uwaterloo.ca}\\[\medskipamount]
{\it \/}                     \>   \\[\bigskipamount]

\end{tabbing}
\end{document}